\documentclass[oneside,11pt]{article}
\usepackage{latexsym}
\usepackage{caption,subcaption}
\usepackage[italian,english]{babel}
\usepackage{cite}

\usepackage{verbatim}
\usepackage{graphicx}
\usepackage{syntonly}
\usepackage{amsfonts}
\usepackage{amssymb}
\usepackage{amsthm}
\usepackage{mathrsfs}
\usepackage[all]{xy}

\def\opn#1#2{\def#1{\operatorname{#2}}}    
\opn\chara{char} \opn\length{\ell}
\opn\projdim{proj\,dim} \opn\injdim{inj\,dim} \opn\rank{rank}
\opn\depth{depth} \opn\grade{grade} \opn\height{height}
\opn\embdim{emb\,dim} \opn\codim{codim}

\opn\Tr{Tr} \opn\bigrank{big\,rank}
\opn\superheight{superheight}\opn\lcm{lcm} \opn\trdeg{tr\,deg}
\opn\reg{reg} \opn\hdeg{hdeg}
\opn\div{div} \opn\Div{Div} \opn\cl{cl} \opn\Cl{Cl}
\opn\Spec{Spec} \opn\Supp{Supp} \opn\supp{supp} \opn\Sing{Sing}
\opn\Ass{Ass}
\opn\Ann{Ann} \opn\Rad{Rad} \opn\Soc{Soc} \opn\Shad{Shad}
\opn\rk{rank} \opn\Ker{Ker} \opn\Coker{Coker} \opn\Im{Im}
\opn\Hom{Hom} \opn\Tor{Tor} \opn\Ext{Ext} \opn\End{End}
\opn\Aut{Aut} \opn\id{id}

\opn\nat{nat}
\opn\pff{pf}  
\opn\Pf{Pf} \opn\GL{GL} \opn\SL{SL} \opn\mod{mod} \opn\ord{ord}
\opn\Poin{Poin}
\opn\aff{aff} \opn\relint{relint} \opn\st{st} \opn\con{conv}
\opn\lk{lk} \opn\cn{cn} \opn\core{core} \opn\vol{vol}
\renewcommand{\con}{\mathrm{conv}}
\opn\gr{gr}

\def\pot#1#2{#1[\kern-0.28ex[#2]\kern-0.28ex]}

\opn\dirlim{\underrightarrow{\lim}}
\opn\invlim{\underleftarrow{\lim}}

\def\Implies{\ifmmode\Longrightarrow \else
     \unskip${}\Longrightarrow{}$\ignorespaces\fi}
\def\implies{\ifmmode\Rightarrow \else
     \unskip${}\Rightarrow{}$\ignorespaces\fi}
\def\iff{\ifmmode\Longleftrightarrow \else
     \unskip${}\Longleftrightarrow{}$\ignorespaces\fi}
\let\:=\colon
\theoremstyle{plain}
\newtheorem{theorem}{Theorem}[section]
\newtheorem{corollary}[theorem]{Corollary}

\newtheorem{proposition}[theorem]{Proposition}
\theoremstyle{definition}
\newtheorem{definition}[theorem]{Definition}
\newtheorem{example}[theorem]{Example}
\newtheorem{remark}[theorem]{Remark}
\newcommand{\eqnsection}{
\renewcommand{\theequation}{{\thesection.\arabic{equation}}}
\makeatletter \csname @addtoreset\endcsname{equation}{section}
\makeatother} \eqnsection
\renewcommand{\Im}{\mathrm{Im}}
\renewcommand{\rk}{\mathrm{rk}}
\renewcommand{\div}{\mathrm{div}}
\renewcommand{\Div}{\mathrm{Div}}
\renewcommand{\Ass}{\mathrm{Ass}}
\renewcommand{\Ann}{\mathrm{Ann}}
\renewcommand{\mod}{\mathrm{mod}}
\renewcommand{\depth}{\mathrm{depth}}

\begin{document}
\date{}

\title{\textbf{\large ON THE BETTI NUMBERS AND\\[1mm] GRACEFULNESS OF SOME PLANAR GRAPHS\\[4mm]}}
\author{\textsc{Maurizio Imbesi}\\[2mm]
\small{\em Department of Mathematical and Computer Sciences, Physical and Earth Sciences}\\
\small{\em University of Messina, Viale F. Stagno d'Alcontres 31, I-98166 Messina, Italy}\\
\small{{\tt e-mail address:} maurizio.imbesi@unime.it}\\[6mm]
\textsc{Monica La\,Barbiera}\\[2mm]
\small{\em Department of Electrical, Electronic and Computer Engineering}\\
\small{\em University of Catania, Viale A. Doria 6, I-95125 Catania, Italy}\\
\small{{\tt e-mail address:} monica.labarbiera@unict.it}
}
\maketitle \thispagestyle{empty}
\noindent {\small \textsc{Abstract.} In this article bipartite planar graphs $St_r$ are investigated, $r$ the number of their plane regions. Bounds for the graded Betti numbers and the projective
dimension of the quotient ring associated to such graphs are discussed. We prove
that $r$ is related to algebraic invariants that arise from the projective resolution of the edge ideal of the graph. We also deal with labeling methods for certain graphs and show that graphs $St_r$ admit a graceful numbering of their edges.}\\[6mm]
\smallskip
\section*{Introduction}
We consider theoretical properties and labeling of planar graphs using computational and commutative algebra methods. We refer to \cite{{B-M},{H}} for general information on planar graphs and \cite{{IL1},{IL2},{IL3}} for algebraic and combinatorial themes examined throughout the paper.
In details, we investigate a class of bipartite planar graphs studied
by Doering and Gunston in \cite{D-G}. Our aim is to study the
constrains for Betti numbers, the projective dimension and the gracefulness of the components of
such graphs using their geometry.\\
In Section 1 we give the necessary definitions and introduce the bipartite graphs
$St_r$, namely planar graphs with $r>1$ regions \mbox{which have vertex}
set $V=\{v_1, \ldots, v_{2r+1} \}$ and edge set $E= \{ \{v_1,v_i\} \ | \
2 \leq i \leq r+1\} \cup \{ \{v_i,v_{i+r} \} \ | \ 2 \leq i \leq
r+1\}\, \cup $ $\{ \{v_i,v_{i+r-1}\} \ | \ 3 \leq i \leq r+1\} \cup \{v_2,v_{2r+1} \}$.\\
In Section 2 we study some aspects of bipartite planar graphs using algebraic methods.
We link the number of the regions of $St_r$ to algebraic invariants of its
edge ideal and in particular, to the graded Betti numbers and the
projective dimension. More precisely, we prove that it is possible
to give an explicit formula for the second graded Betti number in
degree $3$ of these planar graphs linked to the number of their
regions and we give upper bounds for the graded Betti numbers and
projective dimension of the edge ideals of graphs $St_r$.\\
In Section 3 we deal with labeling of some graphs, in particular we examine cycle-related graphs that have graceful labeling and show that a bipartite planar graph $St_r, \,r$ positive integer, admits a graceful numbering of its components. A consequence of graph isomorphisms leads to affirm that the class of Jahangir graphs \cite{LJM} is graceful.

\smallskip
\section{Preliminary notions}
Let $G$ be a graph with a finite vertex set $V =\{v_1,\ldots, v_n\}$
and edge set $E$, that consists of pairs $\{v_i,v_j\}$, called edges,
for some $v_i, v_j \in V$.\\
A \textit{$n$-path} is a graph $P_n$ whose vertices can be listed in the order $v_1, v_2,\dots, v_n$ such that the edges are $\{v_i, v_{i+1}\}, \,i=1,2,\dots, n-1$.\\
A \textit{$n$-cycle} is a graph $C_n$ whose $n$ vertices are connected in a closed chain.\\
A graph with $n$ vertices $v_1,\ldots, v_n$ is said \textit{complete}, and denoted $K_n$, if
there exists an edge for all pairs $\{ v_i,v_j \}$ of vertices of it.\\
A graph is called \textit{bipartite} if its
vertex set $V$ can be partitioned into disjoint subsets
$V_{1}=\{x_{1},\ldots,x_{n}\}$ and $V_{2}=\{y_{1},\ldots,y_{m}\}$,
and any edge joins a vertex of $V_{1}$ to a vertex of $V_{2}$. More, a
bipartite graph is said \textit{complete}, and denoted $K_{n,m}$, if all the vertices of
$V_{1}$ are joined to all the vertices of $V_{2}$.\\[1mm]
Let $G$ be a graph  on vertices $v_1,\ldots, v_n$ and $R = K[X_1,
\ldots , X_n]$ be a polynomial ring over a field $K$, with one
variable $X_i$ for each vertex $v_i$.
\begin{definition}
The \textit{edge ideal $I(G)$} associated to a graph $G$ is the
ideal of $R$ generated by monomials of degree two, $X_iX_j$, on the
$X_1, \ldots, X_n$ variables, such that $\{v_i,v_j\} \in E$ for
$1\leq i, j\leq n$:
\[
I(G) = (\{X_iX_j |\{v_i,v_j\} \in E \}).
\]
\end{definition}
\begin{definition}
Let $G$ be a graph on vertex set $V= \{v_1,\ldots, v_n \}$ and
edge set $E=\{f_1,\ldots, f_q \}$. The \textit{line graph} of
a graph $G$, denoted by $L(G)$, has vertex set equal to the edge set
of $G$ and  two vertices of $L(G)$ are adjacent whenever the
corresponding edges of $G$ have one common vertex.
\[
V(L(G))=E=\{f_1,\ldots, f_q \}\]
\[
E(L(G))= \{ (f_i,f_j) \ | \ f_i=\{v_i,v_j\}, \  f_j=\{v_j,v_k\}, \ i
\neq j, \ j \neq k \}
\]
\end{definition}
\noindent The number of edges of $L(G)$ is given by
\[ |E(L(G))|= - |E| + \sum_{i=1}^n\frac{{\rm deg}^2(v_i)}{2},\]
where $\deg(v_i)$ is the number of edges incident with $v_i$ (see \cite{Vill:Real}).
\begin{definition}
A graph $G$ is \textit{planar} if it has an embedding in the plane, called \textit{plane graph},
such that every two edges are incident only in a vertex of $G$.
\end{definition}
\noindent Every planar graph is divided by its edges in regions, the \textit{faces} of the corresponding plane graph.\\
In showing planarity, it is enough to take graphs whose \emph{connectivity} is at least $2$, the so-called $2$-connected graphs.
\begin{remark} From Euler's polyhedron formula:\\
(a) any planar graph $G$ with $n\geqslant 3$ vertices has at most $3n-6$ edges. Moreover, if $G$ has no triangles, its edges are at most $2n-4$.\\
(b) the complete graph $K_5$ and the complete bipartite graph $K_{3,3}$ are nonplanar; in fact, $K_5$ has $10>9=3n-6$ edges, $K_{3,3}$ with 6 vertices and no triangles, has $9>8=2n-4$ edges.
\end{remark}
\begin{theorem} [Kuratowski] \label{t1}
A graph is planar if and only if it has no subgraph homeomorphic to $K_5$ or $K_{3,3}$.
\end{theorem}
\begin{proof} See \cite{H}, Theorem 11.13\,.
\end{proof}
\begin{remark} \label{R1}
Let $G$ be a planar graph. We call $G$ \emph{maximal} if no edge can be added to it without losing planarity.\\
(1) Any maximal plane graph with $n$ vertices has every face to be a triangle and $3n-6$ edges;\\
(2) A plane graph with $n$ vertices in which every face is a cycle of length $4$, has $2n-4$ edges.
\end{remark}
\noindent From now on, let us consider the class of planar graphs $St_r$ as defined in \cite{D-G}.\\[2mm]
Let $St_r$ be the planar graph with $r>1$ regions
on vertex set $V=\{v_1, \ldots,$ $v_{2r+1} \}$, where $v_1$ is called the \emph{hub}, and edge set $E= \{
\{v_1,v_i\} \ | \ 2 \leq i \leq r+1\} \cup \{ \{v_i,v_{i+r} \} \ | \ 2
\leq i \leq r+1\} \cup \{ \{v_i,v_{i+r-1}\} \ | \ 3 \leq i \leq r+1\} \cup \{v_2,v_{2r+1} \}$.
\begin{remark}\label{R2}
$St_r$ is a bipartite planar graph. The vertex set of $St_r$ can be
partitioned into disjoint subsets $V_{1}=\{x_{1},\ldots,x_{n}\}$ and
$V_{2}=\{y_{1},\ldots,y_{m}\}$, where $n=r+1$ and $m=r$, $m+n=2r+1$,
and any edge joins a vertex of
$V_{1}$ with a vertex of $V_{2}$. The edge set can be written:\\
$E=\{ \{x_1,y_i\} \ | \ 1 \leq i \leq m \} \cup \{ \{x_i,y_{i-1}\} \ | \
2 \leq i \leq n \} \cup \{ \{x_i,y_i\} \ | \ 2 \leq i \leq m \} \cup \{x_n,y_{1}\}$.\\
It follows that $St_r$ is bipartite, in addition it is complete only when $r=2$.
\end{remark}
\begin{example}\label{es1}
Let $St_2$ with  $V=\{v_1, v_2, v_3, v_4, v_5 \}$ and \\
$E=\{\{v_1, v_2\}, \{v_1, v_3\}, \{v_2, v_4\}, \{v_3, v_5\}, \{v_3,
v_4\}, \{v_2, v_5\} \}$. \\
$V$ can be partitioned into disjoint subsets
$\{x_1,x_2,x_3 \} \cup \{y_1,y_2 \}$,\\
where $x_1= v_1$, $x_2= v_4$,  $x_3= v_5$,  $y_1= v_2$,  $y_2=
v_3$.\\
Then: $E=\{\{x_1, y_1\}, \{x_1, y_2\}, \{x_2, y_1\}, \{x_3, y_2\},
\{x_2, y_2\}, \{x_3, y_1 \} \}$.
\begin{figure}[h]
\begin{subfigure}
{0.5\textwidth}
\hspace{1.8cm}
   \includegraphics[width=.5 \linewidth]{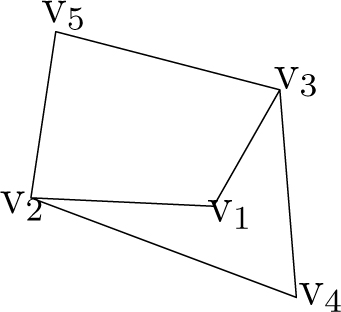}
\label{fig1}
\end{subfigure}
\begin{subfigure}
{0.45\textwidth}
\vspace{6mm}\hspace{1cm}
   \includegraphics[width=.5 \linewidth]{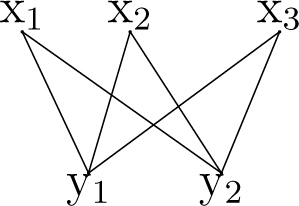}
   \vspace{5mm}
\label{fig2}
\end{subfigure}
\caption{Two representations of the graph $St_2$}
\end{figure}
\end{example}

\section{Syzygy modules of $I(St_r)$}

Let $St_r$ be the  bipartite planar graph with $r>1$ regions. We are
interested to find bounds for the graded Betti numbers that appear
in the minimal graded resolution of its edge ideal. These numbers
determine the rank of the free modules appearing in the minimal
graded resolution  and it is not possible to give a generic formula
to compute them, except for the second and third graded Betti
numbers. But in general, we give upper bounds
for them linked to the number of the regions of  the planar graph $St_r$.\\[1mm]
Let $I(St_r)\subset R$ be the edge ideal of $St_r$. In the following
statement we express the second Betti number of $R/I(St_r)$ in terms of
graph theoretical properties.

\begin{theorem}\label{T1}
Let  $St_r$ be the bipartite planar graph, $r$ be the number of its
regions and $I(St_r)$ be its edge ideal. If
\[
\ldots \rightarrow R^{c}(-4)\oplus R^{b}(-3)\rightarrow R^{q}(-2)
\rightarrow I(St_r) \rightarrow 0
\]
is the minimal graded resolution of $I(St_r)$, then\\
(1) $q= 3r$;\\
(2) $b= \frac{1}{2} r(r+7)$.
\end{theorem}
\begin{proof} $(1)$ $q=|E|= |\{ \{v_1,v_i\} | 2 \leq i \leq r+1\}|
+ |\{ \{v_i,v_{i+r} \} |  2 \leq i \leq r+1\}| + |\{
\{v_i,v_{i+r-1}\} | 3 \leq i \leq r+1\}| + |\{v_2,v_{2r+1} \}|=r
+ r+ (r-1) +1= 3r.$\\
$(2)$ By the formula in \cite{El:Vill},
$b=|E(L(St_r))|- N_3$, where $N_3$ is the number of the triangles of
$St_r$ and $N_3=0$ because the graph is bipartite, one has:\\
$|E(L(St_r))|= - |E| + \sum_{i=1}^N\frac{\deg^2v_i}{2}$, where
$N=2r+1$.\\
Thus, $\sum_{i=1}^{2r+1}\frac{\deg^2v_i}{2}= \frac{r^2}{2} +
r(\frac{3^2}{2}) + r(\frac{2^2}{2}) = \frac{1}{2}(r^2+13r)$,\\
where $\deg v_1=r$, $\deg v_i=3$ for $2\leq i \leq r+1$ and $\deg
v_i=2$ for $r+2 \leq i \leq 2r+1$.\\
Then, $b=|E(L(St_r))|= -3r + \frac{1}{2}(r^2+13r)= \frac{1}{2} r(r+7)$.
\end{proof}
\noindent We study bounds for the graded Betti numbers that appear in the
minimal graded resolution of the edge ideal of $St_r$, in particular
we give upper bounds for them in terms of the number its the
regions.
\begin{definition}
Let $G$ be a graph on vertex set $V(G)$. We call \emph{induced
subgraph} of $G$ the graph $H \subseteq G$ which has an edge between
any two vertices of it if and only if there is an edge between them in $G$.
\end{definition}
\begin{proposition} [\cite{J}, Proposition 4.1.1] \label{3.2}
Let $G$ be a graph. If $H$ is an induced subgraph of $G$ on a subset of
the vertices of $G$, 
then
$$ b_{i_j}(H) \leq  b_{i_j}(G), \quad \forall \ i,j\,,$$
where $b_{i_j}(H)$ (resp. $ b_{i_j}(G)$) are the graded Betti
numbers associated to $H$ (resp. $G$).
\end{proposition}
\begin{proposition}\label{2.2}
Let  $St_r$ be the bipartite planar graph on $2r+1$ vertices, $r$ be
the number of its regions and $\mathcal{I}$ be the edge ideal. Let
$b_{i_j}(St_r)$ be the graded Betti numbers in  the minimal graded
resolution of $R/\mathcal{I}$. Then \[b_{i_j}(St_r) \leq
\sum_{k+l=i+1,k,l\geq 1}{r+1 \choose k}{r+1 \choose l} \    \
\forall \ i  \  \ and  \ \ j=i+1.\]
\end{proposition}
\begin{proof} $St_r$ is a not complete  bipartite planar graph on
two disjoint vertex set $V_1$ and $V_2$, with $|V_1| = r+1$ and
$|V_2| = r$. If follows that $St_r$ is an induced subgraph of the
complete bipartite graph $K_{n,m}$, where $n=m=r+1$, that has a
vertex in more that $St_r$. Then, by Proposition \ref{3.2}
\[
b_{i_j}(St_r) \leq b_{i_j}(K_{n,m}), \    \  \forall  \ i  \  \
\textrm{and}  \ \ j=i+1.\]  By \cite{J}, Theorem 5.2.4,
\[ b_{i_j}(K_{n,m})= \sum_{k+l=i+1, k,l\geq 1}{n \choose k}{m
\choose l}, \  \ \textrm{for} \ \ j=i+1. \]
Hence for $n=m=r+1$ it follows:
\[b_{i_j}(St_r) \leq \sum_{k+l=i+1, k,l\geq 1}{r+1 \choose k}{r+1
\choose l} \    \  \forall  \ i  \  \ \textrm{and}  \ \ j=i+1. \quad \vspace{-1.2cm}
\]
\end{proof}
\medskip
\begin{example}
Let  $St_r$ be the bipartite planar graph on $2r+1$ vertices, $r$ be
the number of its regions and $\mathcal{I}$ be the edge ideal. Let
\[
\ldots \rightarrow  \ldots \oplus R^{d}(-4)\rightarrow
R^{c}(-4)\oplus R^{b}(-3) \rightarrow R^{q}(-2) \rightarrow
\mathcal{I} \rightarrow 0
\]
be the minimal graded resolution of $\mathcal{I}$.
The upper bound for the third Betti number $d$ that appears in the conjecture is:\\[1mm]
$d=b_{3_4}(St_r)\leq \sum_{k+l=4}{r+1 \choose k}{r+1 \choose l}=$\vspace{1mm}\par
$={r+1 \choose 1}{r+1 \choose 3} + {r+1 \choose 2}{r+1 \choose 2}+
{r+1 \choose 3}{r+1 \choose 1}= \frac{1}{12} r(r+1)(7r-4)$.
\end{example}
\noindent Now we give bounds for the projective dimension of the edge ideal of $St_r$.
\begin{definition}
Let $G$ be a graph with vertex set $V$. A subset $\mathcal{A}$ of
$V$ is said minimal vertex cover for $G$ if each edge of $G$ is
incident with one vertex in $\mathcal{A}$ and there is no proper
subset  of $\mathcal{A}$ with this property.
\end{definition}
\begin{definition}
The smallest number of vertices in any minimal vertex cover of $G$
is said vertex covering number. We denote it $\alpha_0(G)$.
\end{definition}
\begin{proposition} [\cite{Vill:Real}] \label{6.1.18}
Let $G$ be a graph and $\mathcal{I}$ be the edge ideal. Then
$\alpha_0(G)= ht(\mathcal{I})$.
\end{proposition}
\begin{theorem}
Let $St_r$ be the bipartite planar graph with $r>1$ regions and
$\mathcal{I}$ be the edge ideal. Then $r-1 \leq
pd_R(\mathcal{I})\leq 2r$.
\end{theorem}
\begin{proof} For the lower bound it is $pd_R(\mathcal{I}) \geq
ht(\mathcal{I})- 1$ (see \cite{Vill:Real}), hence by Proposition
\ref{6.1.18}, $pd_R(\mathcal{I}) \geq \alpha_0(St_r)- 1$.\\
$St_r$ has vertex set $V=\{v_1, \ldots, v_{2r+1} \}$ and edge set
$E= \{ \{v_1,v_i\} | 2 \leq i \leq r+1\} \cup  \{ \{v_i,v_{i+r} \} |
2 \leq i \leq
r+1\} \cup \{ \{v_i,v_{i+r-1}\} | 3 \leq i \leq r+1\} \cup \{v_2,v_{2r+1} \}$.\\
 By definition of
$St_r$ and by its geometry in the plane, it follows that  the
vertices of the minimal vertex cover are all the vertices joined to
$v_1$: $\mathcal{A}(St_r) = \{ v_i | 2 \leq i \leq r+1 \}$. Each
edge of $St_r$ is incident in a vertex of $\mathcal{A}(St_r)$ and
this set is minimal as follows by the description of the edge set.
Hence $\alpha_0(St_r)= r$ and $pd_R(\mathcal{I}) \geq r - 1$.\\
For the upper bound we observe that $St_r$ is a bipartite graph on
vertex set $V=V_1 \cup V_2$ with $|V_1| = r+1$, $|V_2| = r$ and it
is an induced subgraph of the bipartite complete graph $K_{n,m}$
with $n=m=r+1$ as in Proposition \ref{2.2}.\\
The projective dimension of a graph is affected by some simple
transformations of the graphs, such as deleting some edges. Then, as
a consequence of \cite{J}, Proposition 4.1.3, one has
$pd_R(\mathcal{I}) \leq pd_R(I(K_{n,m}))$, with $pd_R(I(K_{n,m}))
=n+m-2$ and $n+m-2=2r$. Hence: $pd_R(\mathcal{I}) \leq 2r$.
\end{proof}

\section{Graceful labeling of graphs $St_r$}

A \emph{graph labeling} is an assignment of integers (labels) to the vertices or edges, or both components, of a graph $G=(V,E)$ under certain conditions. Interesting graph labeling methods are those introduced by Rosa \cite{R} and afterwards by Graham and Sloane \cite{G-S}. Rosa called a function $f: V \longrightarrow \{0,1,\dots, q\}$ a \emph{$\beta$-valuation} if $f$ is an injection such that, when each edge $e=uv$ is assigned the difference $|f(u)-f(v)|$, the resulting edge labels are distinct and are all numbers in the set $\{1,2,\dots, q\}$.\\
Subsequently, Golomb \cite{Go} named \emph{graceful} such labeling, and the graphs which admit a graceful labeling are called \emph{graceful graphs}.\\
Notice that most graphs are not graceful, but graphs that have some sort of regularity in their structure are reasonably graceful. We report an excellent book by Gallian \cite{Ga} that collects an updated dynamic list of results on graph labeling, and the references therein.\\
In this section we show that $St_r$ is a graceful graph and study its labeling explicitly, for all $r>1$.
\begin{definition} A graph $G=(V,E)$ is said to be \emph{graceful} when its vertices are labeled with integers that belong to a subset of $\{0,1,\dots, |E|\}$, and this induces a labeling of the edges of $G$ with all distinct integers from 1 to $|E|$.\\
Thus, $G$ is graceful if and only if there exists an injection that induces a bijection from $E$ to all distinct positive integers up to $|E|$.
\end{definition}
\noindent Important cycle-related graceful planar graphs are
\begin{itemize}
\item Complete bipartite $K_{n,m}$,
\item Fans $F_n=P_n+K_1$,
\item Wheels $W_n=C_n+K_1$,
\item Grids $P_n\times P_m$\,.
\end{itemize}
\noindent For instance, a graceful labeling of a wheel graph $W_n$ is expressed in
\begin{theorem}
Let $W_n$ be a wheel, $n\geq 3$. There exists a graceful numbering for $W_n$, defined by the following sequences, where the hub is assigned $0$:\\[2mm]
for $n$ even, \,$\{0; 2n,1,2n-3,3,2n-5,5,2n-7,7,\dots,2\}$,\\
with last assignment always $2$;\\[2mm]
for $n$ odd, \,$\{0; 2n,1,2n-3,3,2n-5,5,2n-7,7,\dots,2n-2\}$,\\
with last assignment always $2n-2$.
\end{theorem}
\begin{proof}
It can be deduced from \cite{H-K}.
\end{proof}
\noindent Moreover, a graceful labeling of a fan graph $F_n$ is expressed in
\begin{theorem}
Let $F_n$ be a fan, $n\geq 2$. There exists a graceful numbering for $F_n$, defined by the following sequences, where the hub is assigned $0$:\\[2mm]
for $n$ even, \,$\{0; 2n+1,1,2n-1,3,2n-3,5,2n-5,7,\dots,n+2,n\}$,\\
with end always $n$;\\[2mm]
for $n$ odd, \,$\{0; 2n+1,1,2n-1,3,2n-3,5,2n-5,7,\dots,n-1,n+1\}$,\\
with end always $n+1$.
\end{theorem}
\begin{proof}
Refer to \cite{F} and remember that $F_{n-1}$ is isomorphic to $W_n$\,.
\end{proof}
\noindent Now we recall the notion of subdivision of a graph.
\begin{definition} Let $G=(V,E)$. An \emph{edge subdivision operation} for $\{u,v\}\in E$ is the deletion of $\{u,v\}$ from $G$ and the addition of two or more edges $\{u,w_1\},\{w_1,w_2\},\dots,\{w_q,v\}$ along, with the new vertices $w_1,\dots, w_q$\,.\\
This operation generates a new graph $G\,'$, where\\
$G\,'=(V\cup \{w_1,\dots, w_q\}, (E\setminus \{u,v\})\cup \{\{u,w_1\},\{w_1,w_2\},\dots,\{w_q,v\}\})$.\\
A graph $H$ which has been derived from $G$ by a sequence of edge subdivision operations is called a \emph{subdivision} of $G$.
\end{definition}
\noindent Examine other cycle-related planar graphs and their gracefulness.
\begin{definition} [shell] Let $C(n,k)$ denote the cycle $C_n$ with $k$ chords sharing a common endpoint. The unique graph $C(n,n-3)$ is called a \emph{shell}.\\
Notice that $C(n,n-3)$ is the same as the fan graph $F_{n-1}$\,.
\end{definition}
\begin{proposition} [\cite{K-G}] Let $C(n,n-3)$ a shell graph with $n\geq 4$. Then $C(n,n-3)$ is a graceful graph. Moreover, the subdivision graph of $C(n,n-3)$ has a graceful labeling.
\end{proposition}
\begin{definition} [gear] A \emph{gear graph} (or \emph{bipartite wheel graph}), denoted by $G_r$, is a graph obtained from the wheel $W_n$ by adding an extra vertex between every pair of adjacent vertices of the $n$-cycle.\\
Thus $G_r$ has $2r+1$ vertices and $3r$ edges.\\
In addition, we call \emph{multigear graph}, and denote it $(sG)_r$, the graph obtained from $W_n$ by operating a subdivision of every edge of the $n$-cycle into $s$ parts.
\end{definition}
\begin{proposition} [\cite{{M-F},{L}}] \label{gracegear} A gear graph $G_r$ has a graceful labeling.\\
In addition, a multilinear graph $(sG)_r$ is also graceful.
\end{proposition}
\begin{definition} [Jahangir] A \emph{Jahangir graph} $J_{n,m}, n\geq 1, m\geq 3$ consists of a cycle $C_{nm}$ together with a hub $v$ adjacent to cyclically labeled vertices $v_1,v_2,\dots, v_m$ in $C_{nm}$ such that $d(v_i,v_{i+1})=n, \,1\leq i \leq m-1$.\\
Thus $J_{n,m}$ has $nm+1$ vertices and $m(n+1)$ edges.
\end{definition}
\noindent For general information and related studies, take a look at \cite{{G-D},{IL2}}.
\begin{remark} Wheel graphs $W_n$ are part of the family of Jahangir graphs, in that $W_n=J_{1,n}$.
More generally, Jahangir graphs are obtainable as subdivision of wheel graphs; specifically, $J_{n,m}$ derives from $W_m$ by operating an edge subdivision, with the insertion of $(n-1)$ vertices in every edge of $C_m$.
\end{remark}
\begin{remark} \label{gear} There exists an isomorphism between gear graphs $G_r$ and Jahangir graphs $J_{2,r}$\,. The same between $(sG)_r$ and $J_{s,r}$.
\end{remark}
\begin{proposition} A Jahangir graph $J_{n,m}$ has a graceful labeling.
\end{proposition}
\begin{proof} Immediately it descends from Proposition \ref{gracegear} and Remark \ref{gear}.
\end{proof}
\begin{corollary} A bipartite planar graph $St_r, \,r>1$, \,is graceful.
\end{corollary}
\begin{proof} It is easy to show that $St_r$ is isomorphic to the gear graph $G_r$\,.
\end{proof}
\noindent The following result shows gracefulness and illustrates a labeling of $St_r$\,.
\begin{theorem}
Let $St_r, \,r\geq 2$, be a bipartite planar (wheel) graph. A graceful numbering of $St_r$ is:\\[1mm]
for $r$ even,\\
$\{0; 3r,1,3(r-1),4,3(r-2),7,3(r-3),\dots, \widehat{\frac{3r}{2}+1},\dots,3r-5,6,3r-2,3,2\}$,\\
where last assignments are always $3,2$ and \, $\widehat{ }$\, indicates missing vertex;\\[1mm]
for $r$ odd,\\
$\{0; 3r,1,3(r-1),7,3(r-2),5,3(r-3),13,3(r-4),11,3(r-5),19,3(r-6),$ $17,\dots,6,3,3r-5,2\}$,
where last assignments are always $3,3r-5,2$\,.
\end{theorem}
\begin{proof} First of all, we call \emph{$q$-vertex} a vertex of degree $q$ of the graph $St_r$.\\
Let the $2r+1$ vertices of $St_r$ be placed in different lines according to their degree, the hub $w$ ($r$-vertex), the $3$-vertices $v_1,\dots,v_r$ in this order, the $2$-vertices $u_1,\dots,u_r$ in this order.\\
Thereby, of the $3r$ edges of $St_r$, $r$ edges join $w$ with $v_1,\dots,v_r$\,; $2(r-1)$ edges join $v_i$ with $u_i$ and $u_{i+1},\, i=1,\dots, r-1$\,; $2$ edges join $v_r$ with $u_1$ and $u_r$\,.\\
We adopt the strategy of assigning $w$ label $0$, $v_1,v_2$ labels $3r,3(r-1)$, resp., $u_1,u_2$ labels $2,1$, resp..\\
Note that $St_2$ was examined in Example \ref{es1}. It is the unique complete bipartite case $K_{3,2}$ in the class of bipartite planar graphs, so $St_2$ has a graceful labeling given precisely by $\{0,6,1,3,2\}$.\\
Moreover, $St_3$ can be labeled starting from the previous assignments of the vertices and completing the numbering of the other ones. We easily obtain a graceful labeling for it when the remaining $3$-vertex is assigned $4$ and last $2$-vertex is assigned $3$, namely $\{0,9,1,6,3,4,2\}$.\\
$St_4$ is isomorphic to the square planar grid $P_3\times P_3$. A simple labeling scheme for proving that $P_n\times P_2$ is graceful, readily extensible to all grids, was exposed in \cite{M}.\\
Now we will show gracefulness of the graph $St_r=(V,E), \,r\geq 4$.\\[1mm]
\emph{Case 1}: \, $r$ is even.\\
Define $f: V \longrightarrow \{0,1,\dots, 3r\}$ as follows:\\[2mm]
$f(w)\,=\,0$\,;\\[1mm]
$f(v_i)\,=\,3(r+1-i), \;\;\; \forall \,i=1,\dots,r$\,;\\[1mm]
$f(u_j)\,= \left\{ \begin{array}{ll}
2, & \quad {\rm for} \;\; j=1\,;\\
3j-5, & \quad \forall \; j=2,\dots,\frac{r}{2}+1\,;\\
3j-2, & \quad \forall \; j=\frac{r}{2}+2,\dots,r\,. \end{array} \right.$\\[4mm]
\emph{Case 2}: \, $r$ is odd.\\
Define $f: V \longrightarrow \{0,1,\dots, 3r\}$ as follows:\\[2mm]
$f(w)\,=\,0$\,;\\[1mm]
$f(v_i)\,= \left\{ \begin{array}{ll}
3(r+1-i), & \quad \forall \; i=1,\dots,r-1\,;\\
3i-5, & \quad {\rm for} \;\; i=r\,. \end{array} \right.$\\[2mm]
$f(u_j)\,= \left\{ \begin{array}{ll}
2, & \quad {\rm for} \;\; j=1\,;\\
1, & \quad {\rm for} \;\; j=2\,;\\
3j-2, & \quad \forall \; j=3,5,7,\dots,r-2\,;\\
3j-7, & \quad \forall \; j=4,6,8,\dots,r-1\,;\\
3, & \quad {\rm for} \;\; j=r\,. \end{array} \right.$\\[2mm]
In both cases, $f$ induces a function $g$ from $E$ to $\{1,2,\dots,3r\}$ such that each edge $e=uv$ in $E$ is assigned the difference $|f(u)-f(v)|$. It results that distinct edges are assigned distinct differences and, because $|E|=3r$, clearly $g$ is a bijection. Therefore we conclude that $St_r$ is a graceful graph with labeling
given by the following sequences:\\[1mm]
when $r$ is even,\\
$\{0; 3r,1,3(r-1),4,3(r-2),7,3(r-3),\dots, \widehat{\frac{3r}{2}+1},\dots,3r-5,6,3r-2,3,2\}$,\\
where last assignments are always $3,2$ and \, $\widehat{ }$ \, indicates missing vertex;\\[1mm]
when $r$ is odd,\\
$\{0; 3r,1,3(r-1),7,3(r-2),5,3(r-3),13,3(r-4),11,3(r-5),19,3(r-6),$ $17,\dots,6,3,3r-5,2\}$,
where last assignments are always $3,3r-5,2$.
\end{proof}
\begin{example} In the following figures, graceful labelings in different representations of bipartite planar graphs are displayed.\\[-5mm]
\begin{figure}[ht]
\hspace{1.1cm}
\includegraphics[width=.9
\textwidth]{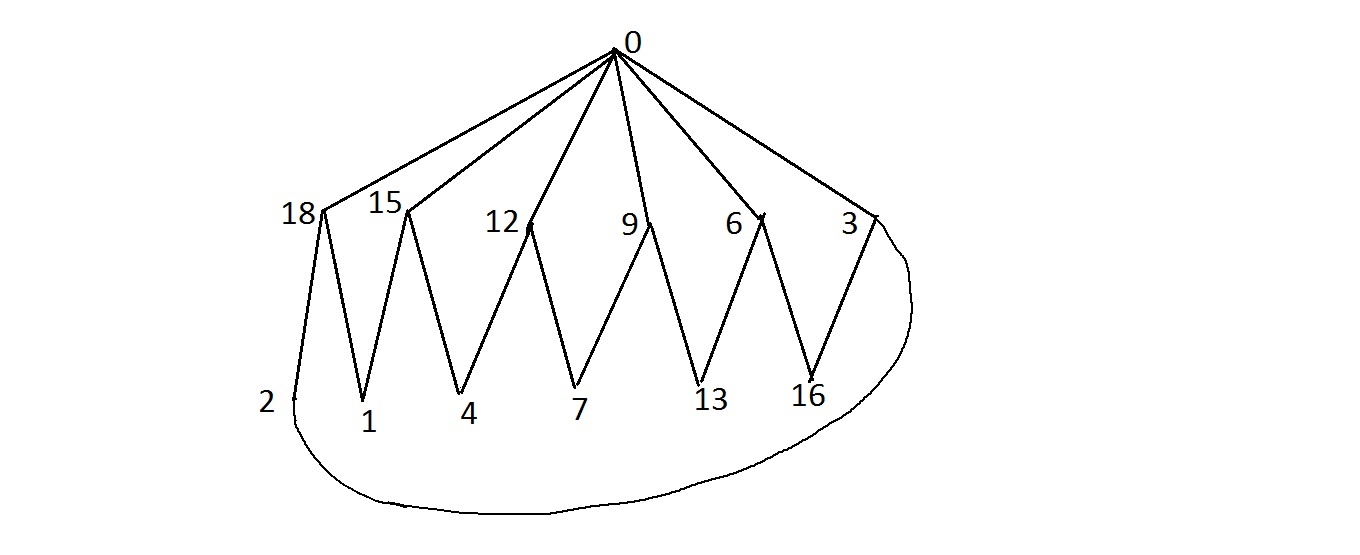}
\vspace{-8mm}
\caption{A graceful labeling of $St_6$ or $G_6$}
\label{St6}
\end{figure}
\vspace{1.1cm}
\begin{figure}[ht]
\hspace{1.6cm}
\includegraphics[width=.9
\textwidth]{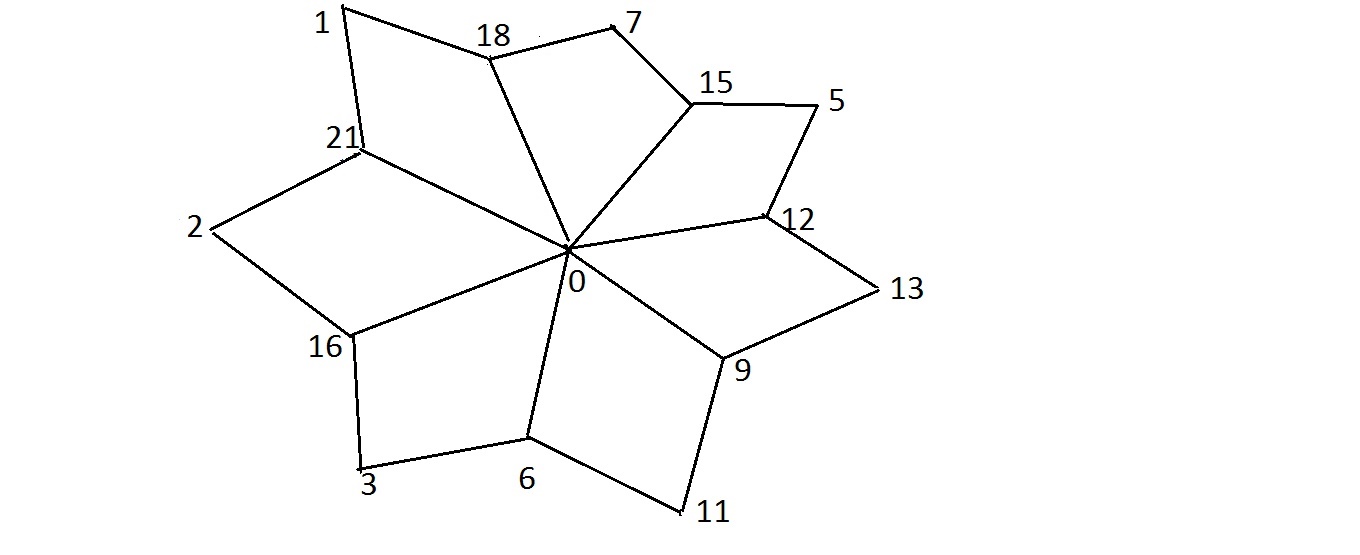}
\vspace{-7mm}
\caption{A graceful labeling of $G_7$ or $J_{2,7}$}
\label{G7}
\end{figure}
\end{example}

\vspace{.5cm} \par\noindent
{\Large \textbf{Acknowledgements}}\\[2.5mm]
The research that led to the present paper was partially supported by a grant of the group GNSAGA of INdAM, Italy.


\end{document}